# Restricted estimation of the cumulative incidence functions corresponding to competing risks

## Hammou El Barmi[1] and Hari Mukerjee[2]

*Baruch College, City University of New York and Wichita State University*

**Abstract:** In the competing risks problem, an important role is played by the cumulative incidence function (CIF), whose value at time $t$ is the probability of failure by time $t$ from a particular type of failure in the presence of other risks. In some cases there are reasons to believe that the CIFs due to various types of failure are linearly ordered. El Barmi et al. [3] studied the estimation and inference procedures under this ordering when there are only two causes of failure. In this paper we extend the results to the case of $k$ CIFs, where $k \geq 3$. Although the analyses are more challenging, we show that most of the results in the 2-sample case carry over to this $k$-sample case.

## 1. Introduction

In the competing risks model, a unit or subject is exposed to several risks at the same time, but the actual failure (or death) is attributed to exactly one cause. Suppose that there are $k \geq 3$ risks and we observe $(T, \delta)$, where $T$ is the time of failure and $\{\delta = j\}$ is the event that the failure was due to cause $j$, $j = 1, 2, \ldots, k$. Let $F$ be the distribution function (DF) of $T$, assumed to be continuous, and let $S = 1 - F$ be its survival function (SF).

The cumulative incidence function (CIF) due to cause $j$ is a sub-distribution function (SDF), defined by

$$(1.1) \qquad F_j(t) = P[T \leq t, \ \delta = j], \ j = 1, 2, \ldots, k,$$

with $F(t) = \sum_j F_j(t)$. The cause specific hazard rate due to cause $j$ is defined by

$$\lambda_j(t) = \lim_{\Delta t \to 0} \frac{1}{\Delta t} P[t \leq T < t + \Delta t, \ \delta = j \mid T \geq t], \ j = 1, 2, \ldots, k,$$

and the overall hazard rate is $\lambda(t) = \sum_j \lambda_j(t)$. The CIF, $F_j(t)$, may be written as

$$(1.2) \qquad F_j(t) = \int_0^t \lambda_j(u) S(u) \, du.$$

Experience and empirical evidence indicate that in some cases the cause specific hazard rates or the CIFs are ordered, i.e.,

$$\lambda_1 \leq \lambda_2 \leq \cdots \leq \lambda_k \ \text{ or } \ F_1 \leq F_2 \leq \cdots \leq F_k.$$

---

[1]Department of Statistics and Computer Information Systems, Baruch College, City University of New York, New York, NY 10010, e-mail: hammou_elbarmi@baruch.cuny.edu

[2]Department of Mathematics and Statistics, Wichita State University, Wichita, KS 67260-0033.

*AMS 2000 subject classifications:* primary 62G05; secondary 60F17, 62G30.

*Keywords and phrases:* competing risks, cumulative incidence functions, estimation, hypothesis test, $k$-sample problems, order restriction, weak convergence.

The hazard rate ordering implies the stochastic ordering of the CIFs, but not *vice versa*. Thus, the stochastic ordering of the CIFs is a milder assumption. El Barmi et al. [3] discussed the motivation for studying the restricted estimation using several real life examples and developed statistical inference procedures under this stochastic ordering, but only for $k = 2$. They also discussed the literature on this subject extensively. They found that there were substantial improvements by using the restricted estimators. In particular, the asymptotic mean squared error (AMSE) is reduced at points where two CIFs cross. For two stochastically ordered DFs with (small) independent samples, Rojo and Ma [17] showed essentially a uniform reduction of MSE when an estimator similar to ours is used in place of the nonparametric maximum likelihood estimator (NPMLE) using simulations. Rojo and Ma [17] also proved that the estimator is better in risk for many loss functions than the NPMLE in the one-sample problem and a simulation study suggests that this result extends to the 2-sample case. The purpose of this paper is to extend the results of El Barmi et al. [3] to the case where $k \geq 3$. The NPMLEs for $k$ continuous DFs or SDFs under stochastic ordering are not known. Hogg [7] proposed a pointwise isotonic estimator that was used by El Barmi and Mukerjee [4] for $k$ stochastically ordered continuous DFs. We use the same estimator for our problem. As far as we are aware, there are no other estimators in the literature for these problems. In Section 2 we describe our estimators and show that they are strongly uniformly consistent. In Section 3 we study the weak convergence of the resulting processes. In Section 4 we show that confidence intervals using the restricted estimators instead of the empiricals could possibly increase the coverage probability. In Section 5 we compare asymptotic bias and mean squared error of the restricted estimators with those of the unrestricted ones, and develop procedures for computing confidence intervals. In Section 6 we provide a test for testing equality of the CIFs against the alternative that they are ordered. In Section 7 we extend our results to the censoring case. Here, the results essentially parallel those in the uncensored case using the Kaplan-Meier [9] estimators for the survival functions instead of the empiricals. In Section 8 we present an example to illustrate our results. We make some concluding remarks in Section 9.

## 2. Estimators and consistency

Suppose that we have $n$ items exposed to $k$ risks and we observe $(T_i, \delta_i)$, the time and cause of failure of the $i$th item, $1 \leq i \leq n$. On the basis of this data, we wish to estimate the CIFs, $F_1, F_2, \ldots, F_k$, defined by (1.1) or (1.2), subject to the order restriction

$$(2.1) \qquad\qquad F_1 \leq F_2 \leq \cdots \leq F_k.$$

It is well known that the NPMLE in the unrestricted case when $k = 2$ is given by (see Peterson, [12])

$$(2.2) \qquad \hat{F}_j(t) = \frac{1}{n} \sum_{i=1}^{n} I(T_i \leq t, \, \delta_i = j), \; j = 1, 2,$$

and this result extends easily to $k > 2$. Unfortunately, these estimators are not guaranteed to satisfy the order constraint (2.1). Thus, it is desirable to have estimators that satisfy this order restriction. Our estimation procedure is as follows.

For each $t$, define the vector $\hat{\mathbf{F}}(t) = (\hat{F}_1(t), \hat{F}_2(t), \ldots, \hat{F}_k(t))^T$ and let $\mathcal{I} = \{\mathbf{x} \in R^k : x_1 \leq x_2 \leq \cdots \leq x_k\}$, a closed, convex cone in $R^k$. Let $E(\mathbf{x}|\mathcal{I})$ denote the least squares projection of $\mathbf{x}$ onto $\mathcal{I}$ with equal weights, and let

$$Av[\hat{\mathbf{F}}; r, s] = \frac{\sum_{j=r}^{s} \hat{F}_j}{s - r + 1}.$$

Our restricted estimator of $F_i$ is

$$(2.3) \qquad \hat{F}_i^* = \max_{r \leq i} \min_{s \geq i} Av[\hat{\mathbf{F}}; r, s] = E((\hat{F}_1, \ldots, \hat{F}_k)^T | \mathcal{I})_i, \quad 1 \leq i \leq k.$$

Note that for each $t$, equation (2.3) defines the isotonic regression of $\{\hat{F}_i(t)\}_{i=1}^k$ with respect to the simple order with equal weights. Robertson et al. [13] has a comprehensive treatment of the properties of isotonic regression. It can be easily verified that the $\hat{F}_i^*$s are CIFs for all $i$, and that $\sum_{i=1}^k \hat{F}_i^*(t) = \hat{F}(t)$, where $\hat{F}$ is the empirical distribution function of $T$, given by $\hat{F}(t) = \sum_{i=1}^n I(T_i \leq t)/n$ for all $t$.

Corollary B, page 42, of Robertson et al. [13] implies that

$$\max_{1 \leq j \leq k} |\hat{F}_j^*(t) - F_j(t)| \leq \max_{1 \leq j \leq k} |\hat{F}_j(t) - F_j(t)| \quad \text{for each } t.$$

Therefore $\parallel \hat{F}_i^* - F_i \parallel \leq \max_{1 \leq j \leq k} \|\hat{F}_j - F_j\|$ for all $i$ where $\|.\|$ is used to denote the sup norm. Since $\|\hat{F}_i - F_i\| \to 0$ $a.s.$ for all $i$, we have

**Theorem 2.1.** $P[\|\hat{F}_i^* - F_i\| \to 0 \ as \ n \to \infty, \ i = 1, 2, \ldots, k] = 1.$

If $k = 2$, the restricted estimators of $F_1$ and $F_2$ are $\hat{F}_1^* = \hat{F}_1 \wedge \hat{F}/2$ and $\hat{F}_2^* = \hat{F}_1 \vee \hat{F}/2$, respectively. Here $\wedge$ ($\vee$) is used to denote max (min). This case has been studied in detail in El Barmi et al. [3].

## 3. Weak convergence

Weak convergence of the process resulting from an estimator similar to (2.3) when estimating two stochastically ordered distributions with independent samples was studied by Rojo [15]. Rojo [16] also studied the same problem using the estimator in (2.3). Praestgaard and Huang [14] derived the weak convergence of the NPMLE. El Barmi et al. [3] studied the weak convergence of two CIFs using (2.3). Here we extend their results to the $k$-sample case. Define

$$Z_{in} = \sqrt{n}[\hat{F}_i - F_i] \quad \text{and} \quad Z_{in}^* = \sqrt{n}[\hat{F}_i^* - F_i], \ i = 1, 2, \ldots, k.$$

It is well known that

$$(3.1) \qquad (Z_{1n}, Z_{2n}, \ldots, Z_{kn})^T \xrightarrow{w} (Z_1, Z_2, \ldots, Z_k)^T,$$

a $k$-variate Gaussian process with the covariance function given by

$$Cov(Z_i(s), Z_j(t)) = F_i(s)[\delta_{ij} - F_j(t)], \quad 1 \leq i, j \leq k, \quad \text{for } s \leq t,$$

where $\delta_{ij}$ is the Kronecker delta. Therefore, $Z_i \stackrel{d}{=} B_i^0(F_i)$ for all $i$, the $B_i^0$s being dependent standard Brownian bridges.

Weak convergence of the starred processes is a direct consequence of this and the continuous mapping theorem. First, we consider the convergence in distribution at a fixed point, $t$. Let

$$(3.2) \qquad \mathcal{S}_{it} = \{j : F_j(t) = F_i(t)\}, \ i = 1, 2, \ldots, k.$$

Note that $\mathcal{S}_{it}$ is an interval of consecutive integers from $\{1, 2, \ldots, k\}$, $F_j(t) - F_i(t) = 0$ for $j \in \mathcal{S}_{it}$, and, as $n \to \infty$,

$$(3.3) \qquad \sqrt{n}[F_j(t) - F_i(t)] \to \infty, \text{ and } \sqrt{n}[F_j(t) - F_i(t)] \to -\infty,$$

for $j > i^*(t)$ and $j < i_*(t)$, respectively, where $i_*(t) = \min\{j : j \in S_{it}\}$ and $i^*(t) = \max\{j : j \in S_{it}\}$.

**Theorem 3.1.** *Assume that* (2.1) *holds and $t$ is fixed. Then*

$$(Z_{1n}^*(t), Z_{2n}^*(t), \ldots, Z_{kn}^*(t))^T \xrightarrow{d} (Z_1^*(t), Z_2^*(t), \ldots, Z_k^*(t))^T,$$

*where*

$$(3.4) \qquad Z_i^*(t) = \max_{i_*(t) \leq r \leq i} \min_{i \leq s \leq i^*(t)} \frac{\sum_{\{r \leq j \leq s\}} Z_j(t)}{s - r + 1}.$$

Except for the order restriction, there are no restrictions on the $F_i$s for the convergence in distribution at a point in Theorem 2. For $k = 2$, if the $F_i$s are *distribution functions* and the $\hat{F}_i$s are the empiricals based on *independent* random samples of sizes $n_1$ and $n_2$, then, using restricted estimators $\hat{F}_i^*$s that are slightly different from those in (2.3), Rojo [15] showed that the weak convergence of $(Z_{1n_1}^*, Z_{2n_2}^*)$ fails if $F_1(b) = F_2(b)$ and $F_1 < F_2$ on $(b, c]$ for some $b < c$ with $0 < F_2(b) < F_2(c) < 1$. El Barmi et al. [3] showed that the same is true for two CIFs. They also showed that, if $F_1 < F_2$ on $(0, b)$ and $F_1 = F_2$ on $[b, \infty)$, with $F_1(b) > 0$, then weak convergence holds, but the limiting process is discontinuous at $b$ with positive probability. Thus, some restrictions are needed for weak convergence of the starred processes.

Let $c_i$ $(d_i)$ be the left (right) endpoint of the support of $F_i$, and let $\mathcal{S}_i = \{j : F_j \equiv F_i\}$ for $i = 1, 2, \ldots, k$. In most applications $c_i \equiv 0$. Letting $i^* = \max\{j : j \in S_i\}$, we assume that, for $i = 1, 2, \ldots, k - 1$,

$$(3.5) \qquad \inf_{c_i + \eta \leq t \leq d_i - \eta} [F_j(t) - F_i(t)] > 0 \text{ for all } \eta > 0 \text{ and } j > i^*.$$

Note that $i \in \mathcal{S}_i$ for all $i$. Assumption (3.5) guarantees that, if $F_j \geq F_i$, then, either $F_j \equiv F_i$ or $F_j(t) > F_i(t)$, except possibly at the endpoints of their supports. This guarantees that the pathology of nonconvergence described in Rojo [15] does not occur. Also, from the results in El Barmi et al. [3] discussed above, if $d_i = d_j$ for some $i \neq j \notin \mathcal{S}_i$, then weak convergence will hold, but the paths will have jumps at $d_i$ with positive probability. We now state these results in the following theorem.

**Theorem 3.2.** *Assume that condition* (2.1) *and assumption* (3.5) *hold. Then*

$$(Z_{1n}^*, Z_{2n}^*, \ldots, Z_{kn}^*)^T \xrightarrow{w} (Z_1^*, Z_2^*, \ldots, Z_k^*)^T,$$

*where*

$$Z_i^* = \max_{i_* \leq r \leq i} \min_{i \leq s \leq i^*} \frac{\sum_{\{r \leq j \leq s\}} Z_j}{s - r + 1}.$$

Note that, if $\mathcal{S}_i = \{i\}$, then $Z_{in}^* \xrightarrow{w} Z_i$ under the conditions of the theorem.

## 4. A stochastic dominance result

In the 2-sample case, El Barmi et al. [3] showed that $|Z_j^*|$ is stochastically dominated by $|Z_j|$ in the sense that

$$P[|Z_j^*(t)| \leq u] > P[|Z_j(t)| \leq u], \ j = 1, 2, \text{ for all } u > 0 \text{ and for all } t,$$

if $0 < F_1(t) = F_2(t) < 1$. This is an extension of Kelly's [10] result for *independent* samples case, but restricted to $k = 2$; Kelly called this result a reduction of stochastic loss by isotonization. Kelly's [10] proof was inductive. For the 2-sample case, El Barmi et al. [3] gave a constructive proof that showed the fact that the stochastic dominance result given above holds even when the order restriction is violated along some contiguous alternatives. We have been unable to provide such a constructive proof for the $k$-sample case; however, we have been able to extend Kelly's [10] result to our (special) dependent case.

**Theorem 4.1.** *Suppose that for some $1 \leq i \leq k$, $S_{it}$, as defined in (3.2), contains more than one element for some $t$ with $0 < F_i(t) < 1$. Then, under the conditions of Theorem 3,*

$$P[|Z_i^*(t)| \leq u] > P[|Z_i(t)| \leq u] \quad \text{for all } u > 0.$$

Without loss of generality, assume that $S_{it} = \{j : F_j(t) = F_i(t)\} = \{1, 2, \ldots, l\}$ for some $2 \leq l \leq k$. Note that $\{Z_i(t)\}$ is a multivariate normal with mean $\mathbf{0}$, and

$$(4.1) \qquad Cov(Z_i(t), Z_j(t)) = F_1(t)[\delta_{ij} - F_1(t)], \quad 1 \leq i, j \leq l.$$

Also note that $\{Z_j^*(t); 1 \leq j \leq k\}$ is the isotonic regression of $\{Z_j(t) : 1 \leq j \leq k\}$ with equal weights from its form in (3.4). Define

$$(4.2) \qquad \mathbf{X}^{(i)}(t) = (Z_1(t) - Z_i(t), Z_2(t) - Z_i(t), \ldots, Z_l(t) - Z_i(t))^T.$$

Kelly [10] shows that, on the set $\{Z_i(t) \neq Z_i^*(t)\}$,

$$(4.3) \qquad P[|Z_i^*(t)| \leq u \,|\, \mathbf{X}^{(i)}] > P[|Z_i(t)| \leq u \,|\, \mathbf{X}^{(i)}] \; a.s. \quad \forall u > 0,$$

using the key result that $\mathbf{X}^{(i)}(t)$ and $Av(\mathbf{Z}(t); 1, k)$ are independent when the $Z_i(t)$'s are independent. Although the $Z_i(t)$s are not independent in our case, they are exchangeable random variables from (4.1). Computing the covariances, it easy to see that $\mathbf{X}^{(i)}(t)$ and $Av(\mathbf{Z}(t); 1, k)$ are independent in our case also. The rest of Kelly's [10] proof consists of showing that the left hand side of (4.3) is of the form $\Phi(a + v) - \Phi(a - v)$, while the right hand side of (4.3) is $\Phi(b + v) - \Phi(b - v)$ using (4.2), where $\Phi$ is the standard normal DF, and $b$ is further away from 0 than $a$. This part of the argument depends only on properties of isotonic regression, and it is identical in our case. This concludes the proof of the theorem.

## 5. Asymptotic bias, MSE, confidence intervals

If $S_{it} = \{i\}$ for some $i$ and $t$, then $Z_i^*(t) = Z_i(t)$ from Theorem 2, and they have the same asymptotic bias and AMSE. If $S_{it}$ has more than one element, then, for $k = 2$, El Barmi et al. [3] computed the exact asymptotic bias and AMSE of $Z_i^*(t)$, $i = 1, 2$, using the representations, $Z_1^* = Z_1 + 0 \wedge (Z_2 - Z_1)/2$ and $Z_2^* = Z_2 - 0 \wedge (Z_2 - Z_1)/2$. The form of $Z_i^*$ in (3.4) makes these computations intractable. However, from Theorem 4, we can conclude that $E[Z_i^*(t)]^2 < E[Z_i(t)]^2$, implying an improvement in AMSE when the restricted estimators are used.

From Theorem 4.1 it is clear that confidence intervals using the restricted estimators will be more conservative than those using the empiricals. Although we believe that the same will be true for confidence bands, we have not been able to prove it.

The confidence bands could always be improved by the following consideration. The $100(1 - \alpha)\%$ simultaneous confidence bands, $[L_i, U_i]$, for $F_i$, $1 \le i \le k$, in the unrestricted case obey the following probability inequality

$$P(F_i \in [L_i, U_i] : 1 \le i \le k) \ge 1 - \alpha.$$

Under our model, $F_1 \le F_2 \le \cdots \le F_k$, this probability is not reduced if we replace $[L_i, U_i]$ by $[L_i^*, U_i^*]$, where

$$L_i^* = \max\{L_j : 1 \le j \le i\} \quad \text{and} \quad U_i^* = \min\{U_j : i \le j \le k\}, 1 \le i \le k.$$

## 6. Hypotheses testing

Let $H_0 : F_1 = F_2 = \cdots = F_k$ and $H_a : F_1 \le F_2 \le \cdots \le F_k$. In this section we propose an asymptotic test of $H_0$ against $H_a - H_0$. This problem has already been considered by El Barmi et al. [3] when $k = 2$, and the test statistic they proposed is $T_n = \sqrt{n} \sup_{x \ge 0} [\hat{F}_2(x) - \hat{F}_1(x)]$. They showed that under $H_0$,

$$(6.1) \qquad \lim_{n \to \infty} P(T_n > t) = 2(1 - \Phi(t)), \quad t \ge 0,$$

where $\Phi$ is the standard normal distribution function.

For $k > 2$, we use an extension of the sequential testing procedure in Hogg [7] for testing equality of distribution functions based on independent random samples. For testing $H_{0j} : F_1 = F_2 = \cdots = F_j$ against $H_{aj} - H_{0j}$, where $H_{aj} : F_1 = F_2 = \cdots = F_{j-1} \le F_j$, $j = 2, 3, \ldots, k$, we use the test statistic $\sup_{x \ge 0} T_{jn}(x)$ where

$$T_{jn} = \sqrt{n} \sqrt{c_j} [\hat{F}_j - Av[\hat{\mathbf{F}}; 1, j - 1]],$$

with $c_j = k(j - 1)/j$. We reject $H_{0j}$ for large values $T_{jn}$, that may be also written as

$$T_{jn} = \sqrt{c_j} [Z_{jn} - Av(\mathbf{Z}_n; 1, j - 1)],$$

where $\mathbf{Z}_n = (Z_{1n}, Z_{2n}, \ldots, Z_{kn})^T$. By the weak convergence result in (3.1) and the continuous mapping theorem, $(T_{2n}, T_{3n}, \ldots, T_{kn})^T$ converges weakly to $(T_2, T_3, \ldots, T_k)^T$, where

$$T_j = \sqrt{c_j} [Z_j - Av[\mathbf{Z}; 1, j - 1]].$$

A calculation of the covariances shows that the $T_j$'s are independent. Also note that

$$T_j \stackrel{d}{=} B_j(F), \quad 2 \le j \le k,$$

where the $B_j$'s are independent standard Brownian motions and $F = \sum_{i=1}^{k} F_i = kF_1$ under $H_0$. We define our test statistic for the overall test of $H_0$ against $H_a - H_0$ by

$$T_n = \max_{2 \le j \le k} \sup_{x \ge 0} T_{jn}(x).$$

By the continuous mapping theorem, $T_n$ converges in distribution to $T$, where

$$T = \max_{2 \le j \le k} \sup_{x \ge 0} T_j(x).$$

Using the distribution of the maximum of a Brownian motion on $[0, 1]$ (Billingsley [2]), and using the independence of the $B_i$'s, the distribution of $T$ is given by

$$P(T \geq t) = 1 - P(\sup_x T_j(x) < t, j = 2, \ldots, k)$$

$$= 1 - \prod_{j=2}^k P(\sup_x B_j(F(x)) < t)$$

$$= 1 - [2\Phi(t) - 1]^{k-1}.$$

This allows us to compute the $p$-value for an asymptotic test.

## 7. Censored case

The case when there is censoring in addition to the competing risks is considered next. It is important that the censoring mechanism, that may be a combination of other competing risks, be independent of the $k$ risks of interest; otherwise, the CIFs cannot be estimated nonparametrically. We now denote the causes of failure as $\delta = 0, 1, 2, \ldots, k$, where $\{\delta = 0\}$ is the event that the observation was censored.

Let $C_i$ denote the censoring time, assumed continuous, for the $i$th subject, and let $L_i = T_i \wedge C_i$. We assume that $C_i$s are identically and independently distributed (IID) with survival function, $S_C$, and are independent of the life distributions, $\{T_i\}$. For the $i$th subject we observe $(L_i, \delta_i)$, the time and cause of the failure. Here the $\{L_i\}$ are IID by assumption.

### 7.1. The estimators and consistency

For $j = 1, 2, \ldots, k$, let $\Lambda_j$ be the cumulative hazard function for risk $j$, and let $\Lambda = \Lambda_1 + \Lambda_2 + \cdots + \Lambda_k$ be the cumulative hazard function of the life time $T$. For the censored case, the unrestricted estimators of the CIFs are the sample equivalents of (1.2) using the Kaplan–Meier [9] estimator, $\widehat{S}$, of $S = 1 - F$:

$$(7.1) \qquad \widehat{F}_j(t) = \int_0^t \widehat{S}(u) \, d\widehat{\Lambda}_j(u), \ j = 1, 2, \ldots, k,$$

with $\widehat{F} = \widehat{F}_1 + \widehat{F}_2 + \cdots + \widehat{F}_k$, where $\widehat{S}$ is chosen to be the left-continuous version for technical reasons, and $\widehat{\Lambda}_j$ is the Nelson–Aalen estimator (see, e.g., Fleming and Harrington, [5]) of $\Lambda_j$. Although our estimators use the Kaplan–Meier estimator of $S$ rather than the empirical, we continue to use the same notation for the various estimators and related entities as in the uncensored case for notational simplicity.

As in the uncensored case, we define our restricted estimator of $F_i$ by

$$(7.2) \qquad \begin{aligned} \widehat{F}_i^* &= \max_{r \leq i} \min_{s \geq i} Av[\widehat{\mathbf{F}}; r, s] \\ &= E((\widehat{F}_1, \ldots, \widehat{F}_k)^T | \mathcal{I})_i, \quad 1 \leq i \leq k. \end{aligned}$$

Let

$$\pi(t) = P[L_i \geq t] = P[T_i \geq t, C_i \geq t] = S(t)S_C(t).$$

Strong uniform consistency of the $\widehat{F}_i^*$s on $[0, b]$ for all $b$ with $\pi(b) > 0$ follows from those of the $\widehat{F}_i$'s [ see, e.g., Shorack and Wellner [18], page 306, and the corrections posted on the website given in the reference] using the same arguments as in the proof of Theorem 2 in the uncensored case.

### 7.2. Weak convergence

Let $Z_{jn} = \sqrt{n}[\hat{F}_j - F_j]$ and $Z_{jn}^* = \sqrt{n}[\hat{F}_j^* - F_j]$, $j = 1, 2, \ldots, k$, be defined as in the uncensored case, except that the unrestricted estimators have been obtained via (7.1). Fix $b$ such that $\pi(b) > 0$. Using a counting process-martingale formulation, Lin [11] derived the following representation of $Z_{in}$ on $[0, b]$:

$$Z_{in}(t) = \sqrt{n} \int_0^t \frac{S(u)dM_i(u)}{Y(u)} - \sqrt{n}F_i(t) \int_0^t \frac{\sum_{j=1}^k dM_j(u)}{Y(u)}$$
$$+ \sqrt{n} \int_0^t \frac{F_i(u)\sum_{j=1}^k dM_j(u)}{Y(u)} + o_p(1),$$

where

$$Y(t) = \sum_{j=1}^n I(L_j \geq t) \text{ and } M_i(t) = \sum_{j=1}^n I(L_j \leq t) - \sum_{j=1}^n \int_0^t I(L_j \geq u)d\Lambda_i(u),$$

the $M_i$'s being independent martingales. Using this representation, El Barmi et al. [3] proved the weak convergence of $(Z_{in}, Z_{2n})^T$ to a mean-zero Gaussian process, $(Z_i, Z_2)$, with the covariances given in that paper. A generalization of their results yields the following theorem.

**Theorem 7.1.** *The process* $(Z_{1n}, Z_{2n}, \ldots, Z_{kn})^T \overset{w}{\Longrightarrow} (Z_1, Z_2, \ldots, Z_k)^T$ *on* $[0, b]^k$, *where* $(Z_1, Z_2, \ldots, Z_k)^T$ *is a mean-zero Gaussian process with the covariance functions, for* $s \leq t$,

$$
\begin{aligned}
(7.3) \quad Cov(Z_i(s), Z_i(t)) &= \int_0^s [1 - F_i(s) - \sum_{j \neq i} F_j(u))][1 - F_i(t) - \sum_{j \neq i} F_j(u))] \frac{d\Lambda_i(u)}{\pi(u)} \\
&+ \sum_{j \neq i} \int_0^s [F_i(u) - F_i(s)][F_i(u) - F_i(t)] \frac{d\Lambda_j(u)}{\pi(u)}.
\end{aligned}
$$

*and, for* $i \neq j$,

$$
\begin{aligned}
(7.4) \quad Cov(Z_i(s), Z_j(t)) &= \int_0^s [1 - F_i(s) - \sum_{l \neq i} F_l(u)][F_j(u) - F_j(t)] \frac{d\Lambda_i(u)}{\pi(u)} \\
&+ \int_0^s [1 - F_j(t) - \sum_{l \neq j} F_l(u)][F_i(u) - F_i(s)] \frac{d\Lambda_j(u)}{\pi(u)} \\
&+ \sum_{l \neq i, j} \int_0^s [F_j(s) - F_j(u)][F_i(t) - F_i(u)] \frac{d\Lambda_l(u)}{\pi(u)}.
\end{aligned}
$$

The proofs of the weak convergence results for the starred processes in Theorems 3.1 and 3.2 use only the weak convergence of the unrestricted processes and isotonization of the estimators; in particular, they do not depend on the distribution of $(Z_1, \ldots, Z_k)^T$. Thus, the proof of the following theorem is essentially identical to that used in proving Theorems 3.1 and 3.2; the only difference is that the domain has been restricted to $[0, b]^k$.

**Theorem 7.2.** *The conclusions of Theorems 3.1 and 3.2 hold for* $(Z_{1n}^*, Z_{2n}^*, \ldots, Z_{kn}^*)^T$ *defined above on* $[0, b]^k$ *under the assumptions of these theorems.*

### 7.3. Asymptotic properties

In the uncensored case, for a $t > 0$ and an $i$ such that $0 < F_i(t) < 1$, if $S_{it} = \{1, \ldots, l\}$ and $l \geq 2$, then it was shown in Theorem 4 that

$$P(|Z_i^*(t)| \leq u) > P(|Z_i(t)| \leq u) \quad \text{for all } u > 0.$$

The proof only required that $\{Z_j(t)\}$ be a multivariate normal and that the random variables, $\{Z_j(t) : j \in S_{it}\}$, be exchangeable, which imply the independence of $\mathbf{X}^{(i)}(t)$ and $Av(\mathbf{Z}(t); 1, l)$, as defined there. Noting that $F_j(t) = F_i(t)$ for all $j \in S_{it}$, the covariance formulas given in Theorem 7.1 show that the multivariate normality and the exchangeability conditions hold for the censored case also. Thus, the conclusions of Theorem 4.1 continue to hold in the censored case.

All comments and conclusions about asymptotic bias and AMSE in the uncensored case continue to hold in the censored case in view of the results above.

### 7.4. Hypothesis test

Consider testing $H_0 : F_1 = F_2 = \cdots = F_k$ against $H_a - H_0$, where $H_a : F_1 \leq F_2 \leq \cdots \leq F_k$, using censored observations. As in the uncensored case, it is natural to reject $H_0$ for large values of $T_n = \max_{2 \leq j \leq k} \sup_{x \geq 0} T_{jn}(x)$, where

$$T_{jn}(x) = \sqrt{n}\sqrt{c_j}[\hat{F}_j(x) - Av(\hat{\mathbf{F}}(x); 1, j-1)]$$
$$= \sqrt{c_j}[Z_{jn}(x) - Av(\mathbf{Z}_n(x); 1, j-1)]$$

with $c_j = k(j-1)/j$, is used to test the sub-hypothesis $H_{0j}$ against $H_{aj} - H_{0j}$, $2 \leq j \leq k$, as in the uncensored case. Using a similar argument as in the uncensored case, under $H_0$, $(T_{2n}, T_{3n}, \ldots, T_{kn})^T$ converges weakly $(T_2, T_3, \ldots, T_k)^T$ on $[0, b]^k$, where the $T_i$'s are independent mean zero Gaussian processes. For $s \leq t$, $Cov(T_i(s), T_i(t))$ simplifies to exactly the same form as in the 2-sample case in El Barmi et al. [3]:

$$Cov(T_i(s), T_i(t)) = \int_0^s S(u)\frac{d\Lambda(u)}{S_C(u)}.$$

The limiting distribution of

$$T_n = \max_{2 \leq j \leq k} \sup_{x \geq 0}[Z_{jn}(x) - Av(\mathbf{Z}_n(x); 1, j-1)]$$

is intractable. As in the 2-sample case, we utilize the strong uniform convergence of the Kaplan–Meier estimator, $\hat{S}_C$, of $S_C$, to define

$$T_{jn}^*(x) = \sqrt{n}\sqrt{c_j}\int_0^t \sqrt{\hat{S}_C(u)}\, d[\hat{F}_j(x) - Av(\hat{\mathbf{F}}(x), 1, j-1)], j = 2, 3, \ldots, k,$$

and define $T_n^* = \max_{2 \leq j \leq k} \sup_{x \geq 0} T_{jn}^*(x)$ to be the test statistic for testing the overall hypothesis of $H_0$ against $H_a - H_0$. By arguments similar to those used in the uncensored case, $(T_{2n}^*, T_{3n}^*, \ldots, T_{kn}^*)^T$ converges weakly to $(T_2^*, T_3^*, \ldots, T_k^*)^T$, a mean zero Gaussian process with independent components with

$$T_j^* \overset{d}{=} B_j(F), \quad 2 \leq j \leq k,$$

where $B_j$ is a standard Brownian motion, and $T_n^*$ converges in distribution to a random variable $T^*$. Since $T_j^*$ here and $T_j$ in the uncensored case (Section 6) have the same distribution, $2 \leq j \leq k$, $T^*$ has the same distribution as $T$ in Section 6, i.e.,

$$P(T^* \geq t) = 1 - [2\Phi(t) - 1]^{k-1}.$$

Thus the testing problem is identical to that in the uncensored case, with $T_n$ of Section 6 changed to $T_n^*$ as defined above. This is the same test developed by Aly et al. [1], but using a different approach.

## 8. Example

We analyze a set of mortality data provided by Dr. H. E. Walburg, Jr. of the Oak Ridge National Laboratory and reported by Hoel [6]. The data were obtained from a laboratory experiment on 82 RFM strain male mice who had received a radiation dose of 300 rads at 5–6 weeks of age, and were kept in a conventional laboratory environment. After autopsy, the causes of death were classified as thymic lymphoma, reticulum cell sarcoma, and other causes. Since mice are known to be highly susceptible to sarcoma when irradiated (Kamisaku et al [8]), we illustrate our procedure for the uncensored case considering "other causes" as cause 2, reticulum cell sarcoma as cause 3, and thymic lymphoma as cause 1, making the assumption that $F_1 \leq F_2 \leq F_3$. The unrestricted estimators are displayed in Figure 1, the restricted estimators are displayed in Figure 2. We also considered the large sample test of $H_0 : F_1 = F_2 = F_3$ against $H_a - H_0$, where $H_a : F_1 \leq F_2 \leq F_3$, using the test described in Section 6. The value of the test statistic is 3.592 corresponding to a $p$-value of 0.00066.

## 9. Conclusion

In this paper we have provided estimators of the CIFs of $k$ competing risks under a stochasting ordering constraint, with and without censoring, thus extending the results for $k = 2$ in El Barmi et al. [3]. We have shown that the estimators are uniformly strongly consistent. The weak convergence of the estimators has been derived. We have shown that asymptotic confidence intervals are more conservative when the restricted estimators are used in place of the empiricals. We conjecture that the same is true for asymptotic confidence bands, although we have not been able to prove it. We have provided asymptotic tests for equality of the CIFs against the ordered alternative. The estimators and the test are illustrated using a set of mortality data reported by Hoel [6].

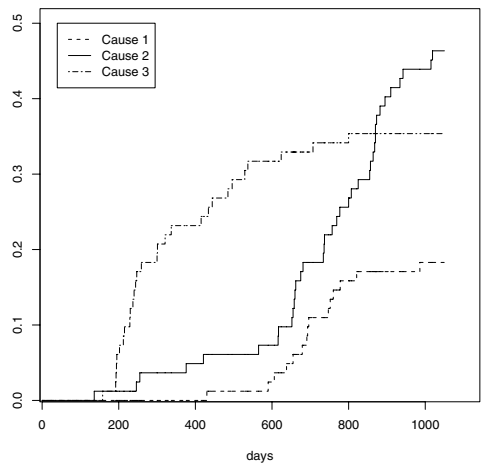

Fig 1. *Unrestricted estimators of the cumulative incidence functions.*

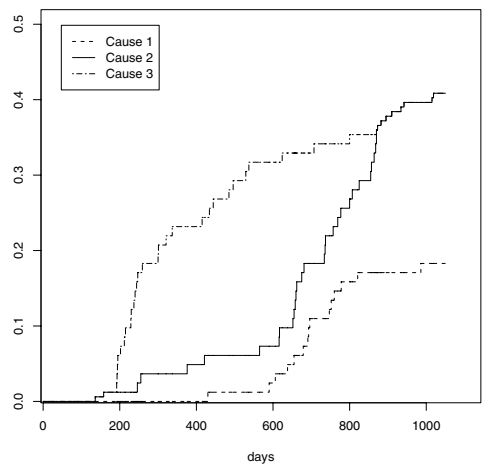

FIG 2. *Restricted estimators of the cumulative incidence functions.*

## Acknowledgments

The authors are grateful to a referee and the Editor for their careful scrutiny and suggestions. It helped remove some inaccuracies and substantially improve the paper. El Barmi thanks the City University of New York for its support through PSC-CUNY.